
\input epsf



\magnification=\magstep1

\documentstyle {amsppt}

\def\a{\alpha}

\def\G{\Gamma}

\def\D{\Delta}

\def\p{\pi}

\def\si{\sigma}
\def\Si{\Sigma}

\def\o{\omega}
\def\O{\Omega}

\def\8{_{\infty}}

\def\tensor{\otimes}

\hfuzz 27pt

\topmatter \title A TQFT for Wormhole
cobordisms \\over the field of rational functions     \endtitle
\rightheadtext{TQFT for Wormholes } \author Patrick M.
Gilmer \endauthor
 \affil Louisiana State University \endaffil
 \address Department of Mathematics,  Baton Rouge, LA
70803 U.S.A  \endaddress
  \email gilmer\@ math.lsu.edu \endemail

  \abstract We consider a cobordism category whose morphisms are punctured
connect sums of $S^1 \times S^2$'s (wormhole spaces) with embedded admissibly
colored  banded trivalent graphs. We define a TQFT on this cobordism category
over  the field of rational functions in an indeterminant $A.$ For  $r$ large,
we recover, by specializing $A$  to a primitive 4rth root of unity, the
Witten-Reshetikhin-Turaev TQFT restricted to links in wormhole spaces.
Thus, for $r$ large, the $r$th Witten-Reshetikhin-Turaev invariant of a link in
some wormhole space,  properly normalized, is the value of a certain rational
function at $e^{\frac{\pi i}{2r}}.$ We relate our work to Hoste  and
Przytycki's calculation of the Kauffman bracket skein module of $S^1 \times
S^2.$ \endabstract
  \endtopmatter

\document
\centerline{7/17/95}
\head Introduction   \endhead

We wish to consider links and graphs in a connected sum of $S^1 \times S^2$'s.
Borrowing a phrase from the physicists, we will use the term `wormhole space'
to  describe this kind of oriented 3-manifold with a relatively simple
topology. Let $\Bbb Q(A)$ denote the field with involution consisting of the
rational functions in the indeterminant $A,$ with involution given by sending
$A$ to $A^{-1}.$

In section one,  we associate to a framed link $L$ in a wormhole space $M$ an
invariant $<L> \in \Bbb Q(A).$
This generalizes the Kauffman Bracket \cite{K} version of the Jones polynomial
\cite{J}. Let $G \subset M $ be a banded trivalent graph in the sense of
\cite{BHMV}\ (4.5)  with an admissible coloring where the set of colors is
taken to be all nonnegative integers.  We will also define $<G> \in \Bbb Q(A).$
Here we make use of Roberts' \cite{R} fusion technique for simplifying bracket
calculations. The well-definedness of the invariant rests  ultimately on the
the well-definedness of the Witten-Reshetikhin-Turaev invariant \cite{W}
\cite{RT}.
See also \cite{L} and \cite{KL}.

In section two, we define the related TQFT using the universal construction
given in \cite{BHMV}. We describe now the cobordism category $\Cal C$ on which
the TQFT is defined.
A punctured wormhole space is a wormhole space with the interiors of some
smooth closed 3-balls deleted. The boundary of a punctured wormhole space is a
collection of 2-spheres. A punctured wormhole space with a properly embedded
admissibly colored graph $(M,G)$ has boundary a disjoint union of 2-spheres
with colored banded points. The objects of $\Cal C$  are disjoint unions
of 2-spheres with colored banded points. We also include the empty suface.  A
morphism from object $\Sigma_1$ to
object $\Sigma_2$ is a triple $(M,G,f).$ Here $f$ is a diffeomorphism of
$\partial (M,G)$ to $-\Sigma_1\coprod \Sigma_2,$ and $M$ is a punctured
wormhole space or a disjoint union of punctured wormhole spaces or is the empty
set.  It
is pleasant that no additional further structure such as a 2-framing or $p_1$
structure is needed. This is because our 3-manifolds are so simple. Since $<\
>$ is a multiplicative involutive invariant, the universal construction
provides a cobordism generated quantization functor from $\Cal C$ to
the category of vector spaces over $\Bbb Q(A).$ We show that the vector spaces
are finite dimensional and that the tensor product axiom holds. Thus we have a
TQFT \cite{A1,A2} which assigns $V(\Sigma)$ to an object $\Sigma$ and $Z_M:
V(\Sigma_1)
\rightarrow V(\Sigma_2)$ to a morphism $M$ from $\Sigma_1$ to $\Sigma_2.$ If
$\Sigma$ is the 2-sphere with $2k$ banded points all colored one, the $\dim(V
(\Sigma) )$ is the $k$th Catalan number.
A version of this TQFT with many fewer objects and morphisms was introduced in
\cite{G,\S 4}.

Hoste and and Przytycki \cite{HP} calculated $\Cal S(S^1 \times S^2),$ the
Kauffman bracket skein module of $S^1 \times S^2.$  $\Cal S(S^1 \times S^2)$ as
a module over $\Bbb Z [A,A^{-1}]$ has a rank one free part, and some torsion.
Thus a framed link in $S^1 \times S^2$ determines a Laurent polynomial in $A.$
In section 3, we show that this polynomial agrees with our $<L>.$  However our
work provides an alternative method of calculating the Hoste-Przytycki
polynomial. This polynomial is also (up to sign) the penultimate
coefficient of $\Gamma(L),$ defined in \cite{G}. The invariant $<L>$ is new for
links in a connected sum of more than one $S^1 \times S^2.$
The fourth section gives two worked examples.

We think that the TQFT introduced here is useful from a pedagogical point of
view
for the following reasons. On the one hand, the objects and morphisms of the
cobordism category are simple topological spaces without any extra structure
such as a 2-framing or $p_1$ structure. On the other hand, the theory is not
completely trivial either.

\head \S 1 Graphs in a wormhole space   \endhead
Let $r$ be an integer greater than two, let $k$ denote $r-2.$  By graph, we
will mean a trivalent banded colored graph as above. Let $G$ be a graph in a
wormhole space $M.$  Suppose that $k$ is greater than all the colors on $G.$
Let $Z_r(M,G)$ denote the Witten-Reshetikhin-Turaev invariant.
This is $\Cal I_{A_r}(M,L)$ evaluated at the element of $\Cal S(S^1 \times
I)^{\tensor \#L}$ given by taking $S_i(\a)$ on a component colored $i$ where
$L$ is the colored link obtained by expanding of $G$ evaluated \cite{L}. We let
$A_r$ denote $e^{\frac{\pi i}{2r}}.$ Also $Z_r(M,G)$ is  $<M,G>_{2p}$ , where
$M$ has been given a $p_1$ structure whose $\si$ invariant is zero
\cite{BHMV}\cite{MV}.
We let $$<G>_r=\frac{Z_{r}(M,G)}{Z_r(M)}$$

$M$ is given by zero framed surgery on the unlink. The number of components of
the unlink is the number of $S^1 \times S^2$ summands.  We mark the components
of the unlink  with dots to indicate  along which components 0-framed 1-surgery
is to be
performed, thinking perhaps of doing a 0-surgery along two three balls on ether
side of the spanning disk \cite{Ki,p.5}. Thus $(M,G)$ maybe
described by  a trivalent colored graph drawn on blackboard with standard
diagram
of unlink whose components have been dotted. We may assume that $G$ is
transverse to the disks which bound the
components of the unlink.
Thus  $Z_r(M,G)$ is the generalized bracket evaluation \cite{KL} of  this
diagram after replacing each component of the unlink with the $\omega_{2r}$ of
\cite{BHMV}. $Z_r(M)$ is simply the generalized
bracket evaluation of the unlink after replacing each dotted component
with $\omega_{2r}.$  Thus if we were to change $\o$ by a scale factor in the
above computation, $<G>_r$ would remain unchanged. Thus we could as well take
$\o$
to be Lickorish's $\o = \sum_i=\D_i S_i(\a)$ or \cite{BHMV}'s $\O_{2r}.$

Suppose the sum of the colors on the strands of $G$ which pass through the
disks is less than or equal to $k.$ Then we may use recoupling theory to
rewrite the  evaluation as a linear combination of evaluations where at most a
single strand passes through each disk. This is called fusion.  Then we may
discard all the terms in the linear combination where the single strand has a
nonzero color using Lickorish's Lemma \cite{L,Lemma 6}. This method of
simplification is due to Roberts \cite{R,Figure 7,Figure 16} \cite{KL,\S12.11}.
In our situation, it is important to see that the above simplification can be
made independent of $r,$ as  long as $k$ is greater than or equal to the sum of
colors passing through each disk. One now has a  linear combination of
trivalent colored graphs
completely unlinked from the unlink. Thus $<G>_r$ is given by the bracket
evaluation of the graph after we delete the unlink. Let $G'$ denote this new
trivalent graph.  $<G>_r$ is given by the bracket evaluation of $G'$ after we
set
$A=A_r.$

Let $<G> \in \Bbb Q(A)$ be the bracket evaluation of $G'.$ Then $<G>_r$ is
simply
$<G>$ evaluated at $A=A_r.$ Since by, say, \cite{BHMV} or \cite{L}, $<G>_r$ is
a well defined isotopy invariant of $G,$ we can conclude that $<G> \in \Bbb
Q(A)$ is also. Here we use the elementary fact that if two rational functions
agree at infinitely many distinct points then they must agree.  This is an
immediate consequence the fundamental theorem of algebra. It is important to
realize that
this argument shows that the result of calculating $<G>$ as above does not
depend on the original surgery description of $M$ or on any choices made in
performing
fusion.

We note that if $L$ is a link diagram, then it describes a framed link in the
wormhole space $S^3.$ If we color this framed link one, we get a graph in $S^3$
which we also denote by $L$. In this situation $<L>$ is the ordinary Kauffman
bracket \cite{K} of the link diagram L, which in turn is version of the Jones
polynomial \cite{J}. Thus $<G> \in \Bbb Q(A)$ generalizes the Kauffman bracket
or Jones polynomial of links in $S^3$ to graphs in a wormhole space.

\proclaim{ Lemma (1.1)} If $G \subset M$ meets an embedded 2-sphere $S \subset
M$ in a single non-zero colored point, then $<G>$ is zero.   \endproclaim

\demo{Proof} Cut $M$ along $S$ and attach two 3-disks to obtain a  space $M'.$
If $M'$ is connected, then $M'$ is a wormhole space and $M$ is obtained  from
$M'$ by a 0-surgery. In this case Lickorish's lemma shows $<G>_r$ to be zero
for $r$ large.
Thus $<G>$ must be zero.

If $M'$ is disconnected, then $M'$ is the disjoint union of two wormhole spaces
and $M$ is obtained by taking their connected sum. Now $<G>_r$ is zero for $r$
large, making use of one basic properties of the Temperley-Lieb  idempotents
$f^{(n)} e_i=0$ \cite{L, Lemma 1}. Thus $<G>$ must be zero.\qed \enddemo

This invariant is extended multiplicatively to disjoint union of wormhole
spaces with graphs.  This invariant is involutive:$<-(M,g)> =
\overline{<(M,g)>}.$

 \head \S 2 A TQFT \endhead

Consider the cobordism category $\Cal C$ described in the introduction. $< \ >$
is an involutive multiplicative invariant on the closed objects of $\Cal C.$
The universal construction of  Blanchet, Habegger,  Masbaum, and  Vogel
then provides a cobordism generated quantization functor $(Z,V)$ to the
category of vector spaces over $\Bbb Q(A).$  To show that $(Z,V)$ is a TQFT we
must show that the tensor product axiom holds, and that the vector spaces
associated to surfaces are finite dimensional
{}.
 Recall elements of $V(\Sigma)$ are equivalence classes of linear combinations
of  graphs in a punctured wormhole space whose boundary is $\Sigma.$  We
 note that $V(\Sigma)$ comes equipped with a nonsingular Hermitian form $<\ ,\
>_\Si$.

\proclaim{ Lemma (2.1)} Let $\Si$ be an object of $\Cal C.$ Let $N$ be a
punctured 3-sphere with boundary the underlying manifold of $\Si.$ $V(\Si)$ is
generated by graphs in $N.$  \endproclaim

\demo{Proof} Let $(N,G)$ be a Wormhole space with graph whose  boundary is
$\Si.$
We may represent $G$ by a graph in a punctured 3-sphere $N'$ together with a
dotted
unlink.  Applying the method of
calculation in \S 1, we obtain a linear combination of graphs  $G'$ in $N'.$
Then
the conjugate linear  form $ <(N,G), \  >_{\Si} =  <(N',G'), \  >_{\Si}.$
Thus  $(N,G)$ and $(N',G')$ represent the same element in $V(\Si).$ \qed
\enddemo

 \proclaim{ Lemma (2.2)} The natural map $V(\Si_1) \tensor V(\Si_2) \rightarrow
 V(\Si_1\coprod\Si_2)$ is surjective. \endproclaim

\demo{Proof} Let $(N,G)$ be a punctured 3-sphere with graph whose  boundary is
$\Si_1\coprod\Si_2.$
Let $S$ be a 2-sphere in $N$ such that $\Si_1$ and $\Si_2$ lie in different
components of $N-S.$ Let $N'$ denote $N$ surgered along $S$.We may use fusion
to represent the same element of  $V(\Si_1\coprod\Si_2)$ as $(N,G)$ by a linear
combination of graphs in $N$ which each meet $S$ in at most a single point. By
Lemma (1.1), we may discard those which meet $S$ in a single point. Let $G'$ be
the resulting linear combination of graphs in $N'.$  $<(N,G), \ >_{\Si} =
<(N',G'), \ >_{\Si}.$
Thus  $(N,G)$ and $(N',G')$ represent the same element in $V(\Si).$  But such
an element is in the image of $V(\Si_1) \tensor V(\Si_2).$\qed
\enddemo

Let $\Sigma$ be a 2-sphere with some banded colored points $\ell.$  $S$ is the
boundary of a 3-ball $B.$
Let $\Cal G$ be an embedded trivalent (noncolored) tree in $B,$ such that $\Cal
G \cap S= \ell.$  Admissible colorings of $G$ give graphs in $B$ and thus
elements of $V(\Si,\ell).$ As in \cite{BHMV} (see also \cite{KL,Chapter 7}),
one may show:

\proclaim{ Lemma (2.3)} The set of admissible colorings of a fixed tree $\Cal
G$ as above forms an orthogonal basis for $V(\Si)$ with respect to the
Hermitian form $<\ ,\ >_\Si$.    \endproclaim

\proclaim{ Theorem (2.4)} $(Z,V)$ is a TQFT.   \endproclaim

\demo{Proof} By Lemmas (2.2) and Lemma (2.3), $V(\Si)$ is finite dimensional
for any
$\Si. $  If we equip $V(\Si_1) \tensor V(\Si_2)$ with the tensor product of the
forms on $V(\Si_1)$ and $V(\Si_2)$, the the map $V(\Si_1) \tensor V(\Si_2)
\rightarrow  V(\Si_1\coprod\Si_2)$ is an isometry , and thus injective. Thus
the map is an isomorphism. \qed\enddemo

\proclaim{Proposition (2.5)} If $\Si$ is a sphere with $2n$ points colored one,
then
$V(\Si)$ is given in the obvious way by the set of diagrams in the disk without
without crossings
and  with boundary $2n$ fixed points. Thus $\dim V(\Si)$ is the nth Catalan
number $c(n)=\frac{1}{n+1}\binom{2n}{n}.$\endproclaim

\demo{Proof} Using the Kauffman relations, it is clear that this set generates
$V(\Si).$ Consider the matrix we get when we pair this set of elements
against itself under $<\ ,\ >_\Si.$ Its determinant is easily seen to be a
polynomial of degree
$n\ c(n)$ in $d = -A^2-A^{-2}.$  So this set of elements is
linearly independent.\qed\enddemo

 \head \S 3 Relation to Hoste and Przytycki's Work \endhead

$\Cal S(S^1 \times S^2)$ modulo its $\Bbb Z[A,A^{-1}]$-torsion submodule is
isomorphic
to $\Bbb Z[A,A^{-1}]$, \cite{HP}.  This quotient can be canonically identified
with $\Bbb Z[A,A^{-1}]. $ The equivalence class of the free generator $1$ given
by the empty link is
identified with one.  Following
Hoste and Przytycki, we let $\p$ denote the projection from $\Cal S(S^1 \times
S^2)$
to this quotient which has canonically identified with $\Bbb Z[A,A^{-1}]. $
Thus
we have $\p:  \Cal S(S^1 \times S^2)\rightarrow \Bbb Z[A,A^{-1}]. $ The
application of
$<\ >$  to links $L$  in $S^1 \times S^2$ colored one defines a $\Bbb
Z[A,A^{-1}]$-module homomorphism to $\Bbb Q(A)$ which must vanish on the
torsion
submodule. Let $\emptyset$ denote the empty link in $S^1 \times S^2,$ then
$<\emptyset> = \p(1) =1 \in \Bbb Z[A,A^{-1}]. $
This proves:

\proclaim{Proposition (3.1)} If $L$ is a framed link in $S^1 \times S^2$
colored one,
then $<L> = \p(L).$\endproclaim

The element of $\Cal S(S^1 \times S^2)$ given by $m$ standardly framed
longitudes is denoted $z^m.$  Let $L_m$ denote this link viewed now as a closed
morphism of $\Cal C$.    We have  by \cite{BHMV,(1.2)}, $<L_m>= \text{Trace}(
Id_{S_m}) = \dim V(S_m).$ Using
Proposition (2.5), we have a new  proof of \cite{HP,Corollary 5}

\proclaim{Corollary (3.2)(Hoste-Przytycki)} $\p(z^{2n+1}) =0,$ and
$\p(z^{2n}) $  is the $n$th Catalan number.\endproclaim

We remark that it follows from Proposition (3.1) that, for a link $L$ colored
one in $S^1 \times S^2,$
$<L> $ actually lies in $\Bbb
Z[A,A^{-1}]$. This also follows from the fact shown in \cite{G} that $\G(\Cal
T)$ has coefficients in $Z[A,A^{-1}].$ The first example in \S4 shows this is
not true for a knot colored one in the connected sum of two $S^1 \times S^2$s.

Although a relationship has  been given between \cite{HP} and
\cite{G, \S4}, the lower bounds given for the wrapping number of links in $S^1
\times S^2$ given by these two papers still seem different, but a detailed
comparison has not been made.

\head \S 4 Two Examples   \endhead

$$\epsffile{fig1.ai}$$

Using the fusion identity on the left of Figure 1, we compute $<K> =
\frac{1}{d}$ for the knot in  the connected sum of two $S^1 \times S^2$s
pictured on the left.

Consider the link $\Cal T$ in $ S^2\times I$ pictured on the left of Figure 2 .
 In
\cite{G}, we calculated a matrix for
$Z_{( S^2\times I,\Cal T)}$ with respect to the basis given in Proposition
(2.5).
The trace of this matrix was
$ A^{-12} - A^{-8} - A^{-4}+1 - 2 A^4 + A^{12} - A^{16}.$  This must be $<L>$
where $L$ is the link in $S^1 \times S^2$  pictured on the right of Figure 2.
$$\epsffile{fig2.ai}$$
 We will now calculate $<L>$ directly by the method described in \S1. We will
make use of the fusion formula given in Figure 3. All unlabelled and undotted
strands are colored one.
$$\epsffile{fig3.ai}$$
Thus $<L> = \frac{1}{d^2}<G_1> + \frac{1}{\D_2}<G_2>,$ where $G_1$ and $G_2$
are the graphs in $S^3$ pictured in Figure 4.
$$\epsffile{fig4.ai}$$
We expand each one using Figure 5.
$$\epsffile{fig5.ai}$$
$<G_1>$ is easily seen to be $d^2(A^2 d+2 + A^{-2}d).$ After we expand $G_2$
as in Figure 3, we may discard the term with coefficient $A^{-2}$ using the
vanishing of the idempotent $f^{(2)}$ times a hook. We are left with two simply
linked theta curves each with edges one, one and two. We may replace these with
loops labelled two. See for instance \cite{KL,page 40}. Thus
$<G_2>$ is $A^2 d+2$ times the evaluation of a Hopf link colored two which
is $A^{-16}\sum_{j=\ 0,\ 2, \ 4}A^{j(j+2)}.$
This is a special case of a formula in \cite{KL,p127} and is also
easily found using \cite{KL,9.15}. Putting this together one calculates that
$<L>$ is indeed the previously calculated trace.


\Refs



\widestnumber \key{BHMV}


 \ref \key A1  \by M. F. Atiyah  \yr 1988 \paper Topological quantum field
theories
\jour Inst. Hautes \'Etudes Sci. Publ. Math.\vol 68  \pages
175-186  \endref

\ref \key A2  \by M. F. Atiyah  \book The geometry and physics of knots, {\rm
Lezioni
Lincee [Lincei Lectures]}\publaddr Cambridge \publ Cambridge Univ. Press
\yr 1990 \endref



 \ref \key  BHMV\by C. Blanchet, N. Habegger, G. Masbaum,  P. Vogel
\paper Topological Quantum Field Theories Derived from the Kauffman Bracket
\jour Topology \toappear \endref

 \ref \key G \by P.~Gilmer
\paper
Invariants for 1-dimensional cohomology
classes arising from TQFT
\paperinfo  preprint \endref




 \ref \key HP \by J.Hoste, J. Przytycki \paper The Kauffman Bracket Skein
Module of  $S^1 \times S^2$ \finalinfo preprint, Preliminary version January
25,1992
 \endref

\ref \key J
\by V.~F.~R. Jones.
\paper A polynomial invariant for knots via von Neumann algebras
\jour Bulletin AMS \vol 12 \pages 103--111 \yr 1985 \endref


 \ref  \key K \bysame \paper State models and the Jones polynomial
\jour Topology\vol 26 \yr 1987 \pages 395-407 \endref



 \ref  \key KL \by L.H. Kauffman, S. Lins \paper Temperley Lieb
recoupling theory and invariants of  $3$-manifolds \jour Annals of Math Studies
 \publaddr Princeton N.J
\publ Princeton Univ. Press
\yr 1994  \endref

\ref \key Ki
\by R.~Kirby
\paper The Topology of 4-manifolds
\jour Springer Lecture Notes in Math.\vol 1374   \endref



 \ref \key  L1 \by W.B.R. Lickorish \paper The skein Method for three manifolds
\jour  J. of Knot Th. and its Ramif. \vol 2 \yr 1993 \pages 171-194  \endref


 \ref \key  MV1 \by  G. Masbaum, P. Vogel
\paper Verlinde Formulae for surfaces with spin structure \inbook Geometric
Topology, Joint U.S. Israel Workshop on Geometric Topology June 10-16, 1992
Technion, Haifa, Israel  {\rm \ Contemporary Mathematics 164} \ed Gordon,
Moriah,
Waynryb \publ American Math. Soc. \yr 1994 \pages 119-137  \endref


 \ref \key RT  \by N. Reshetikhin,  V. Turaev \yr 1991 \paper Invariants of
3-manifolds via link-polynomials and quantum groups \jour Invent.
Math.\vol 103  \pages 547-597  \endref

\ref \key R \by J.Roberts \paper Skein Theory and Turaev-Viro invariants
\jour Topology \toappear \endref


 \ref \key W \by E. Witten\paper Quantum field theory and the Jones
polynomial\jour Commun. Math. Phys.\vol 121 \yr 1989\pages 351-399  \endref



\endRefs
\enddocument
\vfill\eject \end